\documentclass[12pt,twoside]{article}
\usepackage[mathscr]{euscript}
\usepackage{secdot}
\usepackage[utf8]{inputenc}
\usepackage{amsmath, amsthm, amscd, amsfonts, amssymb, color}
\usepackage[bookmarksnumbered,plainpages]{hyperref}

\date{}

\begin{document}

\centerline{}
\centerline {\Large{\bf Characterization of Various Types of Fuzzy Delta Compactness}}
\centerline {\Large{\bf Using Fuzzy Upper Limit of Fuzzy Delta Closed Sets}}
\centerline{}
\centerline{\textbf{Sanjay Roy${^*}$}}
\centerline{Department of Mathematics, Uluberia College}
\centerline{Uluberia, Howrah- 711315,  West Bengal, India}
\centerline{e-mail: sanjaypuremath@gmail.com}
\centerline{}
\centerline{\textbf{Srabani Mondal}}
\centerline{Department of Mathematics, Uluberia College}
\centerline{Uluberia, Howrah- 711315,  West Bengal, India}
\centerline{e-mail: srabanimondal531@gmail.com}
\centerline{}
\centerline{\textbf{T. K. Samanta}}
\centerline{Department of Mathematics, Uluberia College}
\centerline{Uluberia, Howrah- 711315,  West Bengal, India}
\centerline{e-mail: mumpu$_{-}$tapas5@yahoo.co.in}
\centerline{}

\newcommand{\mvec}[1]{\mbox{\bfseries\itshape #1}}
\centerline{}
\newtheorem{Theorem}{\quad Theorem}[section]

\newtheorem{Def}[Theorem]{\quad Definition}

\newtheorem{Thm}[Theorem]{\quad Theorem}

\newtheorem{remark}[Theorem]{\quad Remark}

\newtheorem{corollary}[Theorem]{\quad Corollary}

\newtheorem{note}[Theorem]{\quad Note}

\newtheorem{lemma}[Theorem]{\quad Lemma}

\newtheorem{example}[Theorem]{\quad Example}
\newtheorem{notation}[Theorem]{\quad Notation}

\newtheorem{result}[Theorem]{\quad Result}
\newtheorem{conclusion}[Theorem]{\quad Conclusion}

\newtheorem{proposition}[Theorem]{\quad Proposition}
\newtheorem{prop}[Theorem]{\quad Property}

\begin{abstract}
\textbf{\emph{In this paper, we introduce the concept of various types fuzzy delta $(\delta)$ compactness such as Quasi fuzzy delta compact, Quasi fuzzy countably delta compact, Weakly fuzzy delta compact, $a$-delta compact, Strong fuzzy delta compact, Ultra fuzzy delta compact and  Fuzzy delta compact  and characterize these types of fuzzy delta compactness using the notion of fuzzy upper limit of net of some types of delta $(\delta)$ closed sets. }}
\end{abstract}
{\bf Keywords:}  {Quasi fuzzy delta compact, Quasi fuzzy countably delta compact,Weakly fuzzy delta compact spaces, $a$-delta compact spaces, Strong fuzzy delta compact spaces, Ultra fuzzy delta compact.}\\
\textbf{2020 Mathematics Subject Classification:} 03E72, 54E45, 54D30\\ 
\\${^*}$Corresponding Author


\section{Introduction}
\label{IMV01.S.1}

 In 1965, Zadeh \cite{Zadeh} initiated a novel concept known as fuzzy set as a new mathematical tool for dealing with uncertainties. A fuzzy set $A$ in $X$ is a function from $X$ to $[0, 1]$, that is, each point of $X$  associates to a real number of the interval $[0, 1]$. A fuzzy set $A$ is contained in a fuzzy set $B$  if and only if $A(x)\leq B(x)$ for every $x\in X$ and it is denoted by $A\leq B$.
In last few years many works \cite{Chang, Ganster, Jeon} have been done on fuzzy topological spaces. The concept of fuzzy topology was introduced by Chang \cite{Chang} in 1968. According to Chang,  A fuzzy topology is a family $T$ of fuzzy sets in $X$ which satisfies the following conditions:\\
 $(a)$ $\overline{0},\,\,\, \overline{1}\in T$, where the constant fuzzy set with value $\alpha$ is denoted by $\overline{\alpha},$\\
$(b)$ If $A,\, B\in T$, $A\wedge B\in T$, where $(A\wedge B)(x)= \min\{A(x),\, B(x)\}$,\\
$(c)$ If $A_i \in T$ for each $i\in I$, $\vee_{i\in I} A_i \in T,$ where $(A\vee B)(x)= \max\{A(x),\, B(x)\}$. \\
The pair $(X, \,T)$ called a fuzzy topological
space.  Every member of $T$ is called a open fuzzy set. A fuzzy set is closed if and only if its complement is open. 
 In this paper, we have taken Ming's \cite{Ming} definition of fuzzy point as follows: 
 
 A fuzzy set in  $X$ is called a fuzzy point if it takes the value $0$ for all $y \in X $ except one, say $x \in X $. If its value at $x$ is $\alpha(0 < \alpha\leq 1)$, we denote this fuzzy point by $p_x^\alpha$ or simply by $p$. The fuzzy point $p_x^\alpha$ is said to be contained in a fuzzy set $A$ or to belong to $A$, denoted by $p_x^\alpha\in A$, if $\alpha\leq A(x)$. 

A fuzzy point $p_x^\alpha$ is said to be quasi-coincident \cite{Ming} with $A$ denoted by $p_x^\alpha q A$ if and only if $\alpha + A(x)> 1$. A fuzzy set $A$ is said to be quasi-coincident with $B$, denoted $AqB$, if and only if there exists $x\in X$ such that $A(x)+B(x)> 1$. If $A$ is not quasi-coincident with $B$, then we write $Aq'B$.

A fuzzy set $A$ in a fuzzy topological space $(X, T)$ is called a Q-neighborhood or quasi neighbourhood \cite{Ming} of $p_x^\alpha$ if and only if there exists $B\in T$ such that $p_x^\alpha q B$ and $B\leq A$.

 The notion of compactness plays an important role on a topological space. So, many researchers \cite{Gantner, Georgiou2, Lowen, Lowen2} have analysed the concept of compactness in the field of fuzzy topological spaces. There have the notion of  Quasi fuzzy compact \cite{Georgiou2}, Quasi fuzzy countably compact\cite{Georgiou2, Wong}, Weakly fuzzy compact \cite{Georgiou2}, a-compact \cite{Gantner}, Strong fuzzy compact \cite{Lowen2} and Ultra fuzzy compact \cite{Lowen2}. Then D.N. Georgiou \cite{Georgiou2} has characterized these notions with the help of limit supremum of net of fuzzy closed sets. 
 
 In 1992, M. K. Singal \cite{Singal} introduced a notion of fuzzy delta  open sets and closed sets. With the help of these concepts, we have generalized the notion of fuzzy compact space in this paper and called it fuzzy delta compact space. Also we have introduced the notion of Quasi fuzzy delta compact, Quasi fuzzy countably delta compact, Weakly fuzzy delta compact, a-delta compact, Strong fuzzy delta compact and Ultra fuzzy delta compact and find out the equivalent conditions of these various types of fuzzy delta compactness using fuzzy upper limit of some types of fuzzy net.


\section{Preliminaries }

\begin{Def}\cite{Ming}
Let $(D,\geq)$ be a directed set. Let $X$ be an ordinary set. Let W be the collection of all the fuzzy points in $X$. The function $S:D \rightarrow W $ is called a fuzzy net in $X$. 
\end{Def}

\begin{Def}\cite{Ming}
Let $S=\{S_n:n\in D\}$ be a fuzzy net in $X$. $S$ is said to be Quasi-coincident with $A$ if and only if for each $n\in D$, $S_n$ is quasi-coincident with $A$.
\end{Def}

\begin{Def}\cite{Georgiou}
Let $\{A_n:n\in N\}$ be a net of fuzzy sets in a fuzzy topological space $Y$. Then by $F$-$limsup_N(A_n)$, we denote the fuzzy upper limit of the net $\{A_n:n\in N\}$ in $I^Y$(the collection of all fuzzy sets in $Y$), that is, the fuzzy set which is the union of all fuzzy points $p_x^\alpha$ in $Y$ such that for every $n_o\in N$ and for every fuzzy open $Q$-neighbourhood $U$ of $p_x^\alpha$ in $Y$ there exists an element $n\in N$ for which $n\geq n_o$ and $A_n q U$. In other cases we set $F$-$limsup_N(A_n) = \overline{0}$.

\end{Def}

\begin{Thm}\label{pth1}\cite{Georgiou}
Let $\{A_n:n\in N\}$ be a net of fuzzy closed sets in $Y$ such that $An_1 \leq An_2$ if and only if $n_2 \leq n_1$. Then $F$-$limsup_N(A_n) = \wedge \{A_n:n\in N\}$.
\end{Thm}

\begin{Def}\cite{Noiri}
 An open set U on a topological space $(X,\tau)$is said to be Regular open set if $U=int[cl(U)]$.
\end{Def}

\begin{Def}\cite{Noiri}
 An open set U on a topological space $(X,\tau)$ is said to be Regular closed set if $U=cl[int(U)]$.
\end{Def}

\begin{Def}\cite{Roy}
 A subset U of $(X,\tau)$ is said to be Delta open set if for all $x\in U$, there exists a regular open set P such that $x\in P\subseteq U$.
\end{Def}

\begin{Def}\cite{Roy}
 A subset U of $(X,\tau)$ is said to be Delta Compact set if every delta open cover of U has a finite subcover.
\end{Def}

\begin{Def}\cite{Singal}
A fuzzy subset $A$ of a fuzzy
topological space $(X, T)$ is called a fuzzy regular open set if $int[cl A]=A$. A fuzzy set whose complement is a fuzzy regular open set is called a fuzzy regular closed set. 
\end{Def}

\begin{Def}\cite{Singal}
A subset $A$ of a fuzzy space $(X, T)$ is fuzzy $\delta$-open if for each fuzzy point $p \in A$, there exists a fuzzy regular open set $B$ such that $p\in B\leq A$.\\
A fuzzy set $A$ in $(X, T)$ is fuzzy $\delta$-open iff its
complement is $\delta$-closed.
\end{Def}

\begin{Def}\cite{Georgiou2}
We denote $\sigma$ a class of directed sets. Let $\{A_n:n\in D\}$ be a net of fuzzy sets in $X$. If $D\in \sigma$, then this net is called $\sigma$ net.
\end{Def}

\begin{Def}\cite{Georgiou2}
A fuzzy topological space $X$ is quasi fuzzy compact if  each fuzzy open cover has a finite subcover.
\end{Def}

\begin{Thm}\cite{Georgiou, Georgiou2}
A fuzzy topological space $X$ is  quasi fuzzy  compact if and only if  for every net $\{K_n : n\in D\}$ of fuzzy closed sets in $X$ such that $F$-$limsup_D(K_n) =\overline{0}$, there exists $n_o\in D$ for which $K_n = \overline{0}$ for every $n\in D $, $n\geq n_o$.
\end{Thm}

Let $(X,\tau)$ be a fuzzy topological space and $a\in [0,1)$. Then the family 
     $\{A^{-1}((a,1]):A\in \tau\}$  is a topology on $X$.
This topology is denoted by $i_a(\tau)$.

  			$I_r$ denote the unit interval $I$ equipped with the topology $\delta _r = \{(a,1]:a\in I\}\bigcup\{I\}$. If $\delta$ is a fuzzy topology on $X$, then the initial topology on $X$ for the family of functions $\delta$ and the topological space $I_r$ is denoted by $i(\delta)$.






\section{Quasi Fuzzy Delta Compact Spaces}

\begin{Def}
 A fuzzy topological space X is quasi fuzzy delta compact if each fuzzy delta open cover has a finite subcover.
\end{Def}

\begin{Thm}
A fuzzy topological space $X$ is  Quasi-Fuzzy delta compact if and only if  for every fuzzy net $\{K_n : n\in D\}$ of fuzzy delta closed sets in $X$ with $F$-$limsup_D(K_n) =\overline{0}$, there exists $n_o\in D$ such that $K_n = \overline{0}$ for every $n\in D $, $n\geq n_o$.
\end{Thm}

\textbf{Proof: }
Let $X$ be a quasi fuzzy delta compact space and let $\{K_n:n\in D\}$ be a net of fuzzy delta closed sets in $X$ with $F$-$limsup_D(K_n)= \overline{0}$.
Then for every fuzzy point $p_x^r$ in $X$ there exists a fuzzy open $Q$-neighbourhood $U_{p_x^r}$ of $p_x^{(1-r)}$ in $X$ and an element $n_{p_x^r}\in D$ such that $K_n q' U_{p_x^r}$ for every $n\in D, n\geq n_{p_x^r}$.
So, $K_n(x)+ U_{p_x^r}(x) \leq 1$ for all $x\in X$ and for all $n\geq n_{p_x^r}$,
i.e., $U_{p_x^r}(x)\leq 1-K_n(x)$ for all $x\in X$ and for all $n\geq n_{p_x^r}$,
i.e., $p_x^r \in U_{p_x^r}\leq \overline{1}-K_n$ for all $n\geq n_{p_x^r}$.
Since $K_n$ is a fuzzy delta closed set, $\overline{1} - K_n$ is a fuzzy delta open set for all $n\in D$.
So,there exists a fuzzy regular open set $W_{p_x^r}$ in $X$ such that $p_x^r\in W_{p_x^r}\leq \overline{1}-K_n$  for all $n\geq n_{p_x^r}$.
Clearly, $\{W_{p_x^r} : p_x^r\in X\}$ is a delta open cover of $X$ and
 $\bigvee\{W_{p_x^r}:p_x^r\in X\}=\overline{1}$.
Since the fuzzy space $X$ is quasi fuzzy delta compact, there exist fuzzy points $p_1, p_2,\cdots, p_m \in X$ such that $\bigvee \{W_{p_i} : i= 1, 2,\cdots, m\} = \overline{1}$. 
Let $n_o\in D$ such that $n_o \geq n_{p_i}$ for all $i=1,2,\cdots, n$. Then for every $n\in D, n\geq n_o $, we have $K_n^c\geq \bigvee\{W_{p_i}:i=1, 2,\cdots, m\} = \overline{1}$, 
or, $K_n = \overline{0}$ for every $n\in D,\, n\geq n_o$.\\
Conversely, Let the fuzzy topological space X satisfies the condition of the theorem. 
We prove that the fuzzy topological space X is quasi fuzzy delta compact.
Let $\mathcal{A}$ be a fuzzy delta open cover of X. Let D be the set of all finite subsets of $\mathcal{A}$ directed by inclusion and let $\{K_n : n \in D\}$ be a net of fuzzy delta closed sets in X such that $K_n^c= \bigvee\{A:A\in n\}$.
Obviously, $Kn_{1}\leq Kn_{2}$ if $n_2 \leq n_1$.
Hence $F$-$limsup_D(K_n) = \bigwedge\{K_n:n\in D\}$
                     $= (\bigvee \{K_n^c:n\in D\})^c$
                  $=(\bigvee\{A:A\in \mathcal{A}\})^c$
                  $=\overline{1}^c$
                  $=\overline{0}$.
Thus $F$-$limsup_D(K_n) = \overline{0}$.
By assumption, there exists an element $n_o \in D$ for which $K_n = \overline{0}$ for every $n \in D$, $n \geq n_o $
We have $\overline{1}=K_{n_o}^c=\bigvee\{A:A\in n_o\}$.
Therefore $X$ is quasi fuzzy delta compact.

\begin{Def}
A fuzzy topological space X is called Quasi-Fuzzy $\sigma$ delta compact if for $\sigma$ net $\{K_n : n\in D\}$ of fuzzy delta closed sets in X with $F$-$limsup_D(K_n) =\overline{0}$, there exists $n_o\in D$ such that $K_n = \overline{0}$ for every $n\in D $, $n\geq n_o$.
\end{Def}

\begin{Def}
A fuzzy topological space X is Quasi fuzzy countably delta compact iff each open countable delta cover has a finite subcover. 
\end{Def}

\begin{Thm}
Let $(X,\tau)$ be a fuzzy topological 2nd countable space. The fuzzy topological space X is Quasi fuzzy countably delta compact iff for every net 
$\{K_n: n \in D\}$ of fuzzy delta closed sets in X with $F$-$limsup_D(K_n) = \overline{0}$, there exists 
$n_o \in D$ such that $K_n = \overline{0}$ for every $n\in D$, $n\geqslant n_o$.
\end{Thm}

\textbf{Proof: }
Let $X$ be a quasi fuzzy countably delta compact space and let $\{K_n:n\in D\}$ be a net of fuzzy delta closed sets in $X$ such that $F$-$limsup_D(K_n)= \overline{0}$.
Then for every fuzzy point $p_x^r$ in $X$ there exists a fuzzy open $Q$-neighbourhood $U_{p_x^r}$ of $p_x^{(1-r)}$ in $X$ and an element $n_{p_x^r}\in D$ such that $K_n q' U_{p_x^r}$ for every $n\in D, n\geq n_{p_x^r}$.
So, $K_n(x)+ U_{p_x^r}(x) \leq 1$ for all $x\in X$ and for all $n\geq n_{p_x^r}$,
i.e., $U_{p_x^r}(x)\leq 1-K_n(x)$ for all $x\in X$ and for all $n\geq n_{p_x^r}$,
i.e., $p_x^r \in U_{p_x^r}\leq \overline{1}-K_n$ for all $n\geq n_{p_x^r}$.
Since $K_n$ is a fuzzy delta closed set, $\overline{1} - K_n$ is a fuzzy delta open set for all $n\in D$.
So,there exists a fuzzy regular open set $W_{p_x^r}$ in $X$ such that $p_x^r\in W_{p_x^r}\leq \overline{1}-K_n$  for all $n\geq n_{p_x^r}$.
Since $(X,\tau)$ is 2nd countable, there exists a countable base $\beta$ for $\tau$. So for each $W_{p_x^r}$, there exists $V_{p_x^r}\in \beta$ such that $p_x^r\in V_{p_x^r} \leq W_{p_x^r}\leq \overline{1}-K_n$ for all $n\geq n_{p_x^r}$.
Now $int cl V_{p_x^r}\leq int cl W_{p_x^r} = W_{p_x^r}\leq \overline{1}-K_n$ for all  $n\geq n_{p_x^r}$.
So, $K_nq'intclV_{p_x^r}$ for all $n\geq n_{p_x^r}$.
Clearly, $\{int cl V_{p_x^r} : p_x^r\in X\}$ is a countable delta open cover of $X$ and
 $\bigvee\{int cl V_{p_x^r}:p_x^r\in X\}=\overline{1}$.
Since the fuzzy space $X$ is quasi fuzzy countably delta compact, there exist fuzzy points $p_1, p_2,\cdots, p_m \in X$ such that $\bigvee \{int clV_{p_i}: i= 1,2,\cdots, m\} = \overline{1}$. 
Let $n_o\in D$ such that $n_o \geq n_{p_i}$ for all $i=1,2,\cdots, m$. Then for every $n\in D, n\geq n_o $, we have $K_n^c\geq \bigvee\{int clV_{p_i}: i=1,2,\cdots, m\} = \overline{1}$, 
or, $K_n = \overline{0}$ for every $n\in D,\, n\geq n_o$.\\
Conversely, Let the fuzzy topological space X satisfies the condition of the theorem. 
We prove that the fuzzy topological space X is quasi fuzzy countably delta compact.
Let $\mathcal{A}$ be a countable fuzzy delta open cover of X. Let D be the set of all finite subsets of $\mathcal{A}$ directed by inclusion and let $\{K_n : n \in D\}$ be a net of fuzzy delta closed sets in X such that $K_n^c= \bigvee\{A:A\in n\}$.
Obviously, $Kn_{1}\leq Kn_{2}$ if $n_2 \leq n_1$
Hence $F$-$limsup_D(K_n) = \bigwedge\{K_n:n\in D\}$
                     $= (\bigvee \{K_n^c:n\in D\})^c$
                  $=(\bigvee\{A:A\in \mathcal{A}\})^c$
                  $=\overline{1}^c$
                  $=\overline{0}$.
Thus $F$-$limsup_D(K_n) = \overline{0}$.
By assumption, there exists an element $n_o \in D$ for which $K_n = \overline{0}$ for every $n \in D$, $n \geq n_o$.
We have $\overline{1}=K_{n_o}^c=\bigvee\{A:A\in n_o\}$.
Therefore $X$ is quasi fuzzy countably delta compact.

\begin{Def}\label{d5} A fuzzy topological space $(X,\tau)$ is called Weakly fuzzy delta compact if for every fuzzy delta open cover $\mathcal{U}=\{U_j : j\in J \}$ of fuzzy delta open sets of X, i.e., $\bigvee \{ U_j : j \in J\} = \overline{1}$ and for every $\epsilon > 0$, there exists a finite subfamily $\{Uj_{1},Uj_{2},\cdots,Uj_{m}\}$ of $\mathcal{U}$ such that $\bigvee \{ Uj_{i}:i=1,2,\cdots, m\} \geq \overline{1-\epsilon}$.
\end{Def}

\begin{Thm}
A fuzzy topological space X is weakly delta compact iff for every net $\{K_n : n\in D\}$ of fuzzy delta closed sets in X with $F$-$limsup_D(K_n)= \overline{0}$ and for every $\epsilon> 0$, there exists $n_o \in D$ such that $K_n \leq \overline{\epsilon}$ for every 
$n\in D, n\geq n_o$.
\end{Thm}

\textbf{Proof: }
Let $X$ be a weakly fuzzy delta compact space and let $\{K_n:n\in D\}$ be a net of fuzzy delta closed sets in $X$ such that $F$-$limsup_D(K_n)= \overline{0}$.
Then for every fuzzy point $p_x^r$ in $X$ there exists a fuzzy open $Q$-neighbourhood $U_{p_x^r}$ of $p_x^{(1-r)}$ in $X$ and an element $n_{p_x^r}\in D$ such that $K_n q' U_{p_x^r}$ for every $n\in D, n\geq n_{p_x^r}$.\\
So, $K_n(x)+ U_{p_x^r}(x) \leq 1$ for all $x\in X$ and for all $n\geq n_{p_x^r}$,\\
i.e, $U_{p_x^r}(x)\leq 1-K_n(x)$ for all $x\in X$ and for all $n\geq n_{p_x^r}$,\\
i.e, $p_x^r \in U_{p_x^r}\leq \overline{1}-K_n$ for all $n\geq n_{p_x^r}$.
Since $K_n$ is a fuzzy delta closed set, $\overline{1} - K_n$ is a fuzzy delta open set for all $n\in D$.
So, there exists a fuzzy regular open set $W_{p_x^r}$ in $X$ such that $p_x^r\in W_{p_x^r}\leq \overline{1}-K_n$  for all $n\geq n_{p_x^r}$.
Clearly, $\{W_{p_x^r} : p_x^r\in X\}$ is a delta open cover of $X$ and $\bigvee\{W_{p_x^r}:p_x^r\in X\}=\overline{1}$. 
Since the fuzzy space $X$ is weakly fuzzy delta compact, for any $\epsilon > 0$, there exist fuzzy points $p_1, p_2, \cdots, p_m\in X$ such that $\bigvee \{W_{p_i} : i= 1, 2, \cdots, m\}\geq \overline{1-\epsilon}$. Let $n_o\in D$ such that $n_o \geq n_{p_i}$ for all $i=1, 2, \cdots, m$. Then for every $n\in D, n\geq n_o $, we have $K_n^c\geq \bigvee\{W_{p_i}:i=1, 2,\cdots, m\}\geq \overline{1-\epsilon}$
or, $K_n \leq \overline{\epsilon}$ for every $n\in D,n\geq n_o$.\\
Conversely, Let the fuzzy topological space X satisfy the condition of the theorem. 
We prove that the fuzzy topological space X is weakly fuzzy delta compact.
Let $\mathcal{A}$ be an delta open cover of fuzzy delta open sets of X . Let D be the set of all finite subsets of $\mathcal{A}$ directed by inclusion and let $\{K_n : n \in D\}$ be a net of fuzzy delta closed sets in X such that $K_n^c= \bigvee\{A:A\in n\}$.
Obviously, $Kn_{1}\leq Kn_{2}$ if $n_2 \leq n_1$.
Hence $F$-$limsup_D(K_n) = \bigwedge\{K_n:n\in D\}$
                     $= (\bigvee \{K_n^c:n\in D\})^c$
                  $=(\bigvee\{A:A\in \mathcal{A}\})^c$
                  $=\overline{1}^c$
                  $=\overline{0}$.
Thus $F$-$limsup_D(K_n) = \overline{0}$.
By assumption, there exists an element $n_o \in D$ such that $K_n\leq \overline{\epsilon}$ for every $n \in D$, $n \geq n_o$.
We have $K_{n_o}^c \geq \overline{1-\epsilon}$,
i.e, $\bigvee \{A:A\in n_o\}\geq \overline{1-\epsilon}$.
Therefore $X$ is weakly fuzzy delta compact, as $n_o$ is a finite subset of $\mathcal{A}$.

\begin{Def}
Let $\sigma $ be a class of directed sets. A fuzzy topological space X is called Weakly fuzzy $\sigma $ delta compact if for every $\sigma $ net $\{K_n: n \in D\}$ of fuzzy delta  closed sets in $X$ with $F$-$limsup_D(K_n) = \overline{0}$  and for every $\epsilon > 0$, there exists $n_o \in D $ such that $k_n \leq \overline{\epsilon}$ for every $n \in D,\, n \geq n_o $.
\end{Def}

\begin{Def}
Let $(X,\tau)$be a fuzzy topological space and let $a \in I $. A collection $U$ of fuzzy delta open sets of $X$ will be called an $a$-delta shading of $X$ if for each $x \in X $, there exists $u \in U $ with $u(x)>a$, i.e., $int [cl (u)](x)>a$ or simply $int cl u(x)>a$.
\end{Def}

\begin{Def}
Let $a \in [0,1)$. A fuzzy topological space $(X,\tau)$ is called $a$-delta compact if each $a$-delta shading family in $\tau$ has a finite $a$-delta shading subfamily.
\end{Def}

\begin{Thm}\label{th1}
Let $a \in [0,1)$. A fuzzy topological space $(X, \tau)$ is $a$-delta compact iff for every net $\{ K_n : n \in D\}$ of fuzzy delta closed sets in $X$ with 
$F$-$limsup_D(K_n) <\overline{1-a} $, there exists an element $n_o \in D $ such that $K_n < \overline{1-a}$ for every $n\in D, n \geq n_o$.
\end{Thm}

\textbf{Proof: }
Let $X$ be an $a$-delta compact space and $\{K_n:n\in D\}$ be a net of fuzzy delta closed sets in $X$ such that $F$-$limsup_D(K_n)< \overline{1-a}$. Then for every fuzzy point $p_x^r$ in $X$ with $p_x^r\notin F$-$limsup_D(K_n)$ and  $p_x^r\leq \overline{1-a}$, there exists a fuzzy open $Q$-neighbourhood $Up_x^r$ of $p_x^r$ in $X$ and an element $n_{p_x^r} \in D$ such that $K_n q' U_{p_x^r}$ for every $n\in D, n\geq n_{p_x^r}$, or, equivalently $Up_x^r \leq K_n^c$ for every $n\in D, n\geq n_{p_x^r}$. Since the fuzzy set $Up_x^r$ is a fuzzy open $Q$-neighbourhood of $p_x^r$ in $X$, we have $r+ Up_x^r(x)> 1$, i.e., $Up_x^r(x) > 1-r \geq a$. Now $K_n(x)+ U_{p_x^r}(x) \leq 1$ for all $x\in X$ and for all $n\geq n_{p_x^r}$, or, $U_{p_x^r}(x) \leq (\overline{1}-K_n)(x)$ for all $x\in X$ and for all $n\geq n_{p_x^r}$. i.e., $U_{p_x^r}\leq \overline{1}-K_n$ for all $n\geq n_{p_x^r}$. Let $U_{p_x^r}(x)= \alpha$. Then the fuzzy point $p_x^\alpha \in\overline{1}-K_n$. Since $K_n$ is a fuzzy delta closed set, $\overline{1} - K_n$ is a fuzzy delta open set for all $n\in D$. So, there exists a fuzzy regular open set $W_{p_x^r}$ in $X$ such that $p_x^{\alpha}\in W_{p_x^r}\leq \overline{1}-K_n$ for all $n\geq n_{p_x^r}$. So, $W_{p_x^r}(x)\geq\alpha\geq a$. 
Let $P=\{p_x^r \in X : p_x^r\notin F$-$limsup_D(K_n)$ and $ p_x^r\leq \overline{1-a}\}$. Clearly, the family $\{W_{p_x^r}:p_x^r \in P \}$ is an $a$-delta shading family. Since the fuzzy topological space $X$ is $a$-delta compact, there exists fuzzy points $p_1, p_2, \cdots, p_m \in P$ such that the family $U_1 =\{ W_{p_i}: i=1, 2,\cdots, m\}$ is $a$-delta shading. Clearly, for the family $U_1$ we have $\bigvee \{W_{p_i}:i=1,2,\cdots, m\} > \overline{a}$. Let $n_o\in D$ such that $n_o\geq n_{p_i}$ for every $i=1,2, \cdots, m $. Then for every $n\in D,n\geq n_o$ we have 
$K_n^c \geq \bigvee \{ W_{p_i}:i=1,2,\cdots, m\}> \overline{a}$.
Thus, $K_n < \overline{1-a}$ for every $n\in D, \,n\geq n_o$.\\
Conversely, Let $A'$ be an $a$-delta shading family in $\tau$.
Let D be the set of all finite subsets of $A'$ directed by inclusion.
Let $\{ K_n: n \in D\}$ be a net of fuzzy delta closed sets in $X$ such that $K_n^c = \bigvee \{ A: A\in n\}$ where $n\in D$.\\
Now, $F$-$limsup_D(K_n) = \bigwedge \{ K_n : n\in D\}$
					$= (\bigvee\{K_n^c : n \in D\})^c$
					$= (\bigvee\{A:A \in A'\})^c$.\\		
Since $A'$ is an $a$-delta-shading family, there exists $A \in A'$ such that $A(x)> a$ for all $ x \in X$,
i.e. $int cl A(x) > a$.
Now, $\bigwedge \{K_n : n\in D\}$
			$= (\bigvee\{A:A\in A'\})^c < \overline{1-a}$
and so $F$-$limsup_D(K_n)< \overline{1-a}$. By assumption there exists an element $n_o \in D$	for which $K_n \leq \overline{1-a}$ for every $n\in D , n\geq n_o$.
By the above, we have $K^c_{n_o} =\bigvee \{A:A\in n_o\} > \overline{a}$ 
and therefore the family $\{A:A\in n_o\}$ is an $a$-delta shading subfamily of $A'$. Thus $X$ is $a$-delta compact.

\begin{Def}
Let $a\in [0,1)$. Let $\sigma $ be a class of directed sets. A fuzzy space $X$ is called $\sigma$ $a$-delta compact if for every $\sigma$ net $\{K_n:n\in D\}$ of fuzzy delta closed sets in $X$ with $F$-$limsup_D(K_n) <  \overline{1-a}$, there exists $n_o \in D$ such that $K_n < \overline{1-a}$ for every $n\in D, n\geq n_o$.
\end{Def}

\begin{Thm}\label{th2}
Let $a\in[0,1)$. A fuzzy topological space $(X,\tau)$ is $a$-delta compact iff for every net $\{K_n:n\in D\}$ of delta closed sets of the topological space $(X,i_a(\tau))$ with $limsup_D(K_n)= \phi$
[by $limsup_D(K_n)$ we denote the upper limit of the net $\{K_n:n\in D\}$], there exists an element $n_o\in D$ such that $K_n = \phi $ for every $n\in D, n\geq n_o$.
\end{Thm}

\textbf{Proof: }
Let $(X,\tau)$ be $a$-delta compact space. Let $\{K_n:n \in D\}$ be a net of delta closed sets in $(X,i_a(\tau))$ such that $lim sup _D (K_n) = \phi$. So, for each $x \in X $, there exists a neighbourhood $U_x$ of $x$ and there exists $n_x \in D$ such that $K_n \cap U_x = \phi$ for all $n\geq n_x$ i.e. $U_x\subseteq X-K_n $ for all $n \geq n_x$. Since $X-K_n \in i_a(\tau)$, there exists $A_{n_x} \in \tau$ such that $X-K_n = A_{n_x}^{-1}(a,1]$.
So $U_x\subseteq A_{n_x}^{-1}(a,1]$ for all $n\geq n_x$, or, $A_{n_x}(x) \in (a,1]$ for all $n\geq n_x$.
In particular, $A_{n_x}(x) > a$. So, for each $x \in X$, there exists $A_{n_x} \in \tau$ such that $A_{n_x}(x) > a$. Then the collection $\{A_{n_x}:x\in X\}$ form a $a$-delta shading family. Since $(X,\tau)$ is $a$-delta compact, there exists a finite $a$-delta shading subfamily $\{A_{n_{x_1}}, A_{n_{x_2}}, \cdots, A_{n_{x_k}}\}$. Now there exists $n_o \in D$ such that $n_o\geq n_{x_1},n_{x_2},\cdots,n_{x_k}$. Let $m\geq n_o$ and $x \in X$. Then there exists $n_{x_p} \in \{n_{x_1},n_{x_2},\cdots,n_{x_k}\}$ such that $A_{n_{x_p}} (x)> a$. By our construction, we have $A_{n_{x_p}}(x)>a $ for all $n \geq n_{x_p}$. Since $m\geq n_o\geq n_{x_p}$, $ A_m (x)>a$ for all $m\geq n_o$ and for all $x \in X$.
So, $X\subseteq A_m^{-1}(a,1]$ for all $m\geq n_o$,
i.e., $K_m = \phi $ for all $m\geq n_o$.\\
Conversely, suppose that the fuzzy topological space $X$ satisfies the condition of the theorem. We prove that the fuzzy topological space $X$ is $a$-delta compact.
Let $\mathcal{U} = \{A_\alpha : \alpha \in \Lambda\}$ be $a$-delta shading family in $\tau$. Then for each $x\in X$, there exists $A_{\alpha_x} \in \mathcal{U}$ such that $A_{\alpha_x}(x) >a$ i.e. $int cl A_{\alpha _x} (x) > a$.
We first show that each $A_\alpha^{-1}(a, 1]$ is a delta open in $i_a({\tau})$. Let $x\in A_\alpha^{-1}(a, 1]$. Then $A_\alpha(x)=r>a$. Since $A_\alpha$ is delta open in $\tau$ and $p_x^r\in A_\alpha$, there exists fuzzy regular open set $B_\alpha\in \tau$ such that $p_x^r\in B_\alpha\leq A_\alpha$. Then $A_\alpha(x)\geq B_\alpha(x)\geq r>a$. So, $x\in B_\alpha^{-1}(a, 1]\subseteq A_\alpha^{-1}(a, 1]\cdots (1)$.
Now obviously, $B_\alpha^{-1}(a, 1]\subseteq intcl B_\alpha^{-1}(a, 1]$ as $B_\alpha^{-1}(a, 1]\in i_a(\tau)$.
Let $y \in int cl B_\alpha^{-1}(a,1]$. Then $int cl B_\alpha (y) >a$ which implies that $B_\alpha (y)>a$ as $B_\alpha\in \tau$ is fuzzy regular open set, i.e., $y\in B_\alpha^{-1} (a,1].$ So, $int cl B_\alpha^{-1} (a,1] \subseteq B_\alpha^{-1} (a,1]$. Thus $B_\alpha^{-1}(a,1]$ is regular open in $i_a(\tau)$ for all $\alpha\in \Lambda$. So by $(1)$ $A_\alpha^{-1}(a, 1]$ is a delta open in $i_a(\tau)$.
Let $D$ be the set of all finite subset of $\Lambda$ directed by set inclusion. Then $\cup_{k\in n}A_k^{-1}(a,1]$ is delta open in $i_a(\tau)$ for all $n\in D$.
So $X-\cup_{k\in n}A_k^{-1} (a,1]$ is delta closed in $i_a (\tau)$. Let $K_n = X-\cup_{k\in n}A_{k}^{-1} (a,1]$ for all $n\in D$. Then $\{K_n : n \in D\}$ is a net of delta closed set in $(X,i_a(\tau))$.
Now if $limsup_D K_n \neq \phi$, then for $x \in limsup_D K_n $ there exists $n\in D$ containing  $\alpha_x$ such that $x\in K_n $. So, for every neighbourhood $U$  of $x$, $K_n\cap U\neq\phi$. If we take $U = A_{\alpha_x}^{-1} (a,1]$), the neighbourhood of $x$ then 
$K_n\cap U=(X-\cup_{k\in n}A_{n}^{-1} (a,1])\cap A_{\alpha_x}^{-1} (a,1]$ which is a contradiction as $\alpha_x\in n$.  So, $K_n\cap U=\phi$. So, $limsup_D K_n = \phi$. Then by given condition there exists $n_o\in D$  such that $K_n = \phi $ for all $n \in D$ with $n\geq n_o$. So, $\cup_{k\in n_o}A_{k}^{-1} (a,1] = X $, i.e., for each $x\in X$ there exists $A_k$ where $k\in n_o$ such that $A_k(x)> a$. So we get $a$-delta shading subfamily $\{A_k:\,k\in n\}$ of $\mathcal{U}$ in $\tau$. Hence the fuzzy topological space $X$ is $a$-delta compact.

\begin{Def}
A fuzzy topological space $(x,\delta)$ is called Strong fuzzy delta compact iff it is $a$-delta compact for every $a\in [0,1)$.
\end{Def}

\begin{Thm}
 A fuzzy topological space $X$ is strong fuzzy delta compact iff for every $a\in[0,1)$ and for every net $\{K_n:n\in D \}$ of fuzzy delta closed sets in $X$ with 
 $F$-$limsup_D(K_n)< \overline{1-a}$, there exists an element $n_o \in D $ such that $K_n< \overline{1-a}$ for every $n\in D,n\geq n_o$.
\end{Thm}

\textbf{Proof: } Follows from the Theorem $\ref{th1}$.

\section{ULTRA-FUZZY DELTA COMPACT SPACES}
\begin{Def}
A fuzzy topological space $(X,\tau)$ is called ultra-fuzzy delta compact if $(X,i(\tau))$ is delta compact.
\end{Def}

\begin{Def}
A fuzzy topological space $(X,\tau)$ is called ultra-fuzzy $\sigma$ delta compact if $(X,i(\tau))$ is $\sigma$ delta compact.
\end{Def}

\begin{Thm}
A fuzzy topological space $(X,\tau)$ is ultra-fuzzy delta compact iff for every net $\{K_n:n\in D\}$ of fuzzy delta closed sets of the topological space $(X,i(\tau))$ such that $limsup_D(K_n)= \phi$, there exists an element $n_o\in D$ such that $K_n = \phi$ for every $n\in D, n\geq n_o$.
\end{Thm}
\textbf{Proof: } Follows from the Theorem $\ref{th2}$.
\section{FUZZY DELTA COMPACT SPACES}
\begin{Def}
A fuzzy topological space $(X,\tau)$ is called fuzzy delta compact if for every family $\mathcal{U}$ of fuzzy delta open sets of $X$ and for every $a\in I$ with $\bigvee\{U:U\in \mathcal{U}\} \geq \overline{a}$ and for every $\epsilon \in (0,a]$ there exists a finite subfamily $\mathcal{U}_1$ of $\mathcal{U}$ such that $\bigvee\{U:U\in \mathcal{U}_1\} \geq \overline{a-\epsilon}$.
\end{Def}

\begin{Thm}
Let $X$ be a fuzzy delta compact space. Then for every net $\{K_n:n\in D\}$ of fuzzy delta closed sets in $X$ with $F$-$limsup_D(K_n)=\overline{0}$, there exists $n_o \in D$ such that $K_n = \overline{0}$ for every $n\in D,\,n\geq n_o$.
\end{Thm}

\textbf{Proof: } Let $X$ be a fuzzy delta compact space and let $\{K_n:n\in D\}$ be a net of fuzzy delta closed sets in $X$ such that $F$-$limsup_D(K_n)= \overline{0}$.
Then for every fuzzy point $p_x^\lambda$ in $X$ there exists a fuzzy  open $Q$-neighbourhood $U_{p_x^\lambda}$ of $p_x^{1-\lambda}$ in $X$ and an element $n_{p_x^\lambda} \in D$ such that $K_nq'U_{p_x^\lambda}$ for every $n\in D, n\geq n_{p_x^\lambda}$.
So, $K_n(x)+ U_{p_x^\lambda}(x)\leq 1$ for all $x\in X$ and for all $n\geq n_{p_x^\lambda}$,
or, $U_{p_x^\lambda}(x) \leq (\overline{1}-K_n)(x)$ for all $x\in X$ and for all $n\geq n_{p_x^\lambda}$,
i.e, $p_x^\lambda \in U_{p_x^\lambda}\leq \overline{1}-K_n$ for all $n\geq n_{p_x^\lambda}$.
Since $K_n$ is a fuzzy delta closed set, $\overline{1} - K_n$ is a fuzzy delta open set for all $n\in D$.
So, there exists a fuzzy regular open set $W_{p_x^\lambda}$ in $X$ such that $p_x^\lambda\in W_{p_x^\lambda}\leq \overline{1}-K_n$ for all $n\geq n_{p_x^\lambda}$, that is,  $K_nq'W_{p_x^\lambda}$ for all $n\geq n_{p_x^\lambda}$. Clearly, $\{W_{p_x^\lambda} : p_x^\lambda\in X\}$ is a delta open cover of $X$ and
$\bigvee\{W_{p_x^\lambda}:p_x^\lambda\in X\}=\overline{1}$.
Since the fuzzy space $X$ is fuzzy delta compact, for any $\epsilon >0$, there exists fuzzy points $p_1, p_2, \cdots, p_m \in X$ such that $\bigvee \{W_{p_i} : i= 1, 2, \cdots, m\}\geq\overline{1-\epsilon}$.
Let $n_o\in D$ such that $n_o \geq n_{p_i}$ for all $i=1, 2,\cdots, m$.
Then for every $n\in D,\,n\geq n_o $, we have $K_n q' \bigvee\{W_{p_i}:i=1, 2, \cdots, m\}$,  or, $K_n q' \overline{1-\epsilon}$,  i.e, $K_n(x) \leq \epsilon $ for all $x \in X $ and $n\geq n_o$. Thus $K_n= \overline{0}$ for every $n\in D,\,n\geq n_o$.
 
\begin{Thm}
A fuzzy topological space $X$ is fuzzy delta compact iff for every $a \in I$, for every net $\{K_n:n\in D\}$ of fuzzy delta closed sets in $X$ with $F$-$limsup_D(K_n) \leq \overline{1-a}$ and for every $\epsilon \in (0,a]$, there exists an element $n_o \in D$ such that $K_n \leq \overline{1-a+\epsilon}$ for every $n\in D, \,n\geq n_o$.

\textbf{Proof: }
Let $X$ be a fuzzy delta compact space, $a\in I ,\{K_n:n\in D\}$ be a net of fuzzy delta closed sets in $X$ such that $F$-$limsup_D(K_n) \leq \overline{1-a}$ and $\epsilon\in (0,a]$. Then for every fuzzy point $p_x^r$ in $X$ for which $F$-$limsup_D(K_n) < p_x^r\leq \overline{1-a+(\epsilon/2)}$, there exists a fuzzy open $Q$-neighbourhood $U_{p_x^r}$ of $p_x^r$ in $X$ and an element $n_{p_x^r} \in D$ such that $K_n q' U_{p_x^r}$ for every $n\in D, n\geq n_{p_x^r}$, or, equivalently $U_{p_x^r}\leq K_n^c$ for every $n\in D, n\geq n_{p_x^r}$. Since the fuzzy set $U_{p_x^r}$ is a fuzzy open $Q$-neighbourhood of $p_x^r$ in $X$ we have $r+U_{p_x^r} > 1$, i.e., $U_{p_x^r}(x) > 1-r \geq {a-\epsilon/2}$. Now $K_n(x)+ U_{p_x^r}(x) \leq 1$ for all $x\in X$ and for all $n\geq n_{p_x^r}$,
or, $U_{p_x^r}(x) \leq (\overline{1}-K_n)(x)$ for all $x\in X$ and for all $n\geq n_{p_x^r}$. Let $U_{p_x^r}(x)=\alpha$. Then $p_x^\alpha \in U_{p_x^r}\leq \overline{1}-K_n$ for all $n\geq n_{p_x^r}$. Since $K_n$ is a fuzzy delta closed set, $\overline{1} - K_n$ is a fuzzy delta open set for all $n\in D$.
So, there exists a fuzzy regular open set $W_{p_x^r}$ in $X$ such that $p_x^\alpha\in W_{p_x^r}\leq \overline{1}-K_n$ for all $n\geq n_{p_x^r}$, 
i.e., $K_nq'W_{p_x^r}$ for all $n\geq n_{p_x^r}$.
Let $P= \{p_x^r \in X : F$-$limsup_D(K_n)< p_x^r \leq \overline{1-a+(\epsilon/2)}\}$.
Clearly, $\bigvee \{W_{p_x^r} : p_x^r \in P\}\geq \overline{a-(\epsilon/2)}$. Since the fuzzy space $X$ is fuzzy delta compact and $\epsilon/2\in (0,a-\epsilon/2]$, there exists fuzzy points $p_1, p_2, \cdots, p_m \in X $ such that $\bigvee\{ W_{p_i}:i=1,2,\cdots, m\} \geq \overline{a-\epsilon/2-\epsilon/2}$, or, $\bigvee\{W_{p_i}:i=1,2,\cdots, m\} \geq \overline{a-\epsilon}$.
Let $n_o \in D$ be such that $n_o\geq n_{p_i}$ for every $i=1,2, \cdots, m$. Then for every $n\in D,\,n\geq n_o$ we have $K_n^c \geq \bigvee \{W_{p_i}:i=1,2,\cdots, m\}\geq \overline{a-\epsilon}$. Thus $K_n \leq \overline{1-a+\epsilon}$ for every $n\in D, \,n\geq n_o$.\\
Conversely, suppose that the fuzzy topological space $X$ satisfies the condition of the theorem. We prove that the fuzzy topological space $X$ is fuzzy delta compact.
Let $a\in I$, $\mathcal{A}$ be a family of fuzzy delta open sets such that $\bigvee\{U:U\in \mathcal{A}\}\geq \overline{a}$ and $\epsilon\in (0,a]$. Let $D$ be the set of all finite subsets of $\mathcal{A}$ directed by inclusion and let $\{K_n:n\in D\}$ be a net of fuzzy closed sets in $X$ such that $K_n^c = \bigvee\{A: A\in n\}$ where $n\in D$. Now it is easy to see that $n_1\leq n_2$ iff $n_1\subseteq n_2$ iff $K_{n_1}^c\leq K_{n_2}^c$ iff $K_{n_2}\leq K_{n_1}$. So by the Theorem \ref{pth1}, $F$-$limsup_D(K_n) = \bigwedge\{K_n:n\in D\}$. Also we have $\bigwedge \{K_n:n\in D\} = (\bigvee \{K_n^c:n\in D\})^c=(\bigvee \{A:A\in \mathcal{A}\})^c$
For the family $\mathcal{A}$ we have $\bigvee \{A:A\in \mathcal{A}\}\geq \overline{a}$. Thus $\bigwedge \{K_n:n\in D\} = (\bigvee \{A:A\in \mathcal{A}\})^c \leq \overline{1-a}$.
and therefore $F$-$limsup_D(K_n)\leq \overline{1-a}$. By assumption there exists an element $n_o\in D$ for which  $K_n \leq \overline{1-a+\epsilon}$ for every $n\in D,\, n\geq n_o$. By the above we have $K^c_{n_o}= \bigvee \{A:A\in n_o\}\geq \overline{a-\epsilon}$ and therefore the family $\{A:A\in n_o\}$ is a finite subfamily of $\mathcal{A}$ such that $\bigvee \{A:A\in n_o\}\geq \overline{a-\epsilon}$.
Hence the fuzzy space $X$ is fuzzy delta compact.
\end{Thm}

\textbf{Conclusion: }  In this paper, we have studied various types of fuzzy delta compact spaces and seen that every types of delta compact spaces can be characterized by limit supremum of the net of delta closed sets. We know that in a general topological space, every $\delta$-closed subset of a $\delta$-compact space is $\delta$-compact and  a $\delta$-compact subset of a $T_3$ space is $\delta$-closed set \cite{Roy}. So, one can try to characterize a fuzzy delta closed set by a collection of fuzzy delta compact sets.

\end{document}